\documentclass[reqno,twoside,a4paper,11pt]{amsart}

\usepackage{amsthm, amssymb, amsmath, amsfonts,eufrak,mathrsfs,dsfont,stmaryrd}
\usepackage{url}
\usepackage[colorlinks=true,pdfstartview=FitV,pdftex]{hyperref}
\hypersetup{linkcolor=black,citecolor=blue,urlcolor= blue}
\usepackage{verbatim}
\usepackage{multirow}
\usepackage{color}
\usepackage[utf8]{inputenc}
\usepackage[displaymath, mathlines]{lineno}
\IfFileExists{epsf.def}{\input epsf.def}{\usepackage{epsf}}
\usepackage{mathtools}
\usepackage{graphicx, xcolor}
\usepackage{fancyhdr}
\usepackage{fancyheadings}

\IfFileExists{myowntimes.sty}{%
\usepackage{mathpazo}
\usepackage{mathrsfs}

}
{}
\DeclareFontFamily{OT1}{eusb}{} \DeclareFontShape{OT1}{eusb}{m}{n} {<5> <6> <7> <8> <9> <10> <11> <12> <14.4> eusb10}{}
\DeclareMathAlphabet{\eusb}{OT1}{eusb}{m}{n}

\DeclareFontFamily{OT1}{eusm}{} \DeclareFontShape{OT1}{eusm}{m}{n} {<5> <6> <7> <8> <9> <10> <11> <12> <14.4> eusm10}{}
\DeclareMathAlphabet{\eusm}{OT1}{eusm}{m}{n}

\DeclareFontFamily{OT1}{eufm}{} \DeclareFontShape{OT1}{eufm}{m}{n} {<5> <6> <7> <8> <9> <10> <11> <12> <14.4> eufm10}{}
\DeclareMathAlphabet{\mathfrak}{OT1}{eufm}{m}{n}

\DeclareFontFamily{OT1}{fraktura}{}
\DeclareFontShape{OT1}{fraktura}{m}{n} {<5> <6> <7> <8> <9> <10> <11> <12> <13> <14.4> [1.1] eufm10}{}
\DeclareMathAlphabet{\fraktura}{OT1}{fraktura}{m}{n}

\DeclareFontFamily{OT1}{cmfi}{} \DeclareFontShape{OT1}{cmfi}{m}{n} {<5> <6> <7> <8> <9> <10> <11> <12> <13> <14.4> [0.9] cmfi10}{}
\DeclareMathAlphabet{\cmfi}{OT1}{cmfi}{b}{n}

\DeclareFontFamily{OT1}{cmss}{} \DeclareFontShape{OT1}{cmss}{m}{n} {<5> <6> <7> <8> <9> <10> <11> <12> <13> <14.4> cmss10}{}
\DeclareMathAlphabet{\cmss}{OT1}{cmss}{m}{n}

\DeclareMathAlphabet{\mathpzc}{OT1}{pzc}{m}{it}

\newtheoremstyle{thm}{1.8ex}{1.8ex}{\itshape\rmfamily}{} {\bfseries\rmfamily}{}{2ex}{}

\newtheoremstyle{def}{1.8ex}{1.8ex}{\slshape\rmfamily}{} {\bfseries\rmfamily}{}{2ex}{}

\newtheoremstyle{rem}{1.8ex}{1.8ex}{\rmfamily}{} {\bfseries\rmfamily}{}{2ex}{}

\theoremstyle{thm}

\newcommand\cour[1]{{\fontfamily{pcr}\selectfont #1}}

\theoremstyle{thm}
\newtheorem{theorem}{Theorem}[section]
\newtheorem{lemma}[theorem]{Lemma}
\newtheorem{proposition}[theorem]{Proposition}

\newtheorem*{Main Theorem}{Main Theorem.}

\newtheorem*{special theorem}{Lindeberg-Feller Theorem for Martingales}

\newtheorem{remark}[theorem]{Remark}
\newtheorem{remarks}[theorem]{Remarks}

\theoremstyle{def}

\theoremstyle{rem}

\numberwithin{equation}{section}

\renewcommand{\section}{\secdef\sct\sect}
\newcommand{\sct}[2][default]{%
\refstepcounter{section}
\addcontentsline{toc}{section}{{\tocsection {}{\thesection}{\!\!\!\!#1\dotfill}}{}}
\vspace{0.7cm}
\centerline{\scshape\thesection.\ #1} \nopagebreak \vspace{0.2cm}}
\newcommand{\sect}[1]{%
\vspace{0.4cm} \centerline{\large\scshape\rmfamily #1}
\vspace{0.2cm}}

\renewcommand{\subsection}{\secdef\subsct\sbsect}
\newcommand{\subsct}[2][default]{\refstepcounter{subsection}
\addcontentsline{toc}{subsection}
{{\tocsection{\!\!}{\hspace{1.2em}\thesubsection}{\!\!\!\!#1\dotfill}}{}}
\nopagebreak\vspace{0.45\baselineskip} {\flushleft\bf
\thesubsection~\bf #1.~}
\\*[3mm]\noindent
\nopagebreak}
\newcommand{\sbsect}[1]{\vspace{0.1cm}\noindent
\textbf{#1.~}\vspace{0.1cm}}

\renewcommand{\subsubsection}{%
\secdef \subsubsect\sbsbsect}
\newcommand{\subsubsect}[2][default]{%
\refstepcounter{subsubsection} 
\addcontentsline{toc}{subsubsection}{{\tocsection{\!\!}
{\hspace{3.05em}\thesubsubsection}{\!\!\!\!#1\dotfill}}{}}
\nopagebreak
\vspace{0.15\baselineskip} \nopagebreak {\flushleft\rmfamily
\itshape\thesubsubsection
\ \rmfamily #1\/.}\ }
\newcommand{\sbsbsect}[1]{\vspace{0.1cm}\noindent
\rmfamily \itshape
\arabic{section}.\arabic{subsection}.\arabic{subsubsection} \
\sffamily #1\/.\ }

\renewcommand{\caption}[1]{%
\vglue0.5cm
\refstepcounter{figure}
\begin{minipage}{0.9\textwidth}\small {\sc Figure~\thefigure. }#1\end{minipage}}

\def\myffrac#1#2 in #3{\raise 2.6pt\hbox{$#3 #1$}\mkern-1.5mu\raise 0.8pt\hbox{$#3/$}\mkern-1.1mu\lower 1.5pt\hbox{$#3 #2$}}

\definecolor{lightgray}{gray}{0.5}

\newcommand{\twocite}[2]{\cite{#1}--\cite{#2}}
\newcommand{\treecite}[3]{\cite{#1}--\cite{#2}--\cite{#3}}

\newcommand{\N}        {\mathbb N}

\newcommand{\R}        {\mathbb R}

\newcommand{\PP}        {\mathbb P}

\newcommand{\Z}        {\mathbb Z}

\newcommand{\E}        {\mathbb E}

\newcommand{\scrC}      {\mathscr{C}} 

\newcommand{\scrH}      {\mathscr{H}}

\newcommand{\twoeqref}[2]{\eqref{#1}--\eqref{#2}}
\newcommand{\treeeqref}[3]{\eqref{#1}--\eqref{#2}--\eqref{#3}}

\newcommand{\diam}      {\hbox{{\rm diam}}}

\newcommand{\1}{\mathds 1}

\title[Return Probability and Traps]
{
Anomalous Heat Kernel for Random Walks 
\\
in Random Environments of Conductances
}
\author[O. Boukhadra]
{
Omar~Boukhadra
}


\begin{document} 

\maketitle
\vspace{-0,3cm}
\centerline{Department of Mathematics, University of Constantine 1}
\medskip
\centerline{\small\cour{\url{boukhadra@umc.edu.dz}}}

\begin{abstract} 
\textit{
We study the trapping phenomenon of random walks in random environments of i.i.d. random conductances on the bonds of the grid $\mathbb{Z}^d$, the so-called random conductance model. Our main results concern the important model with conductances in $[0, 1]$ and a polynomial-tailed law near zero for which we find the correct order of decay of the anomalous heat kernel for $d \ge 5$. In $d = 4$, the behavior is found to be normal. In addition, we look at the symmetrical situation with conductances in $[1, \infty)$ with a polynomial law at infinity, which also shows opposite return probability behaviors.
}
\end{abstract}

\hypersetup{linkcolor=red}

\section{\bf Introduction}
\label{}

Let us consider the Random Conductance Model (RCM) on the grid $\Z^d$ where $o$ denotes the origin. 
First, set $E_d = \{ e = \{x, y\}: x, y \in \Z^d, |x - y| =1 \}$ where $| \cdot |$ is the Euclidean distance in $\R^d$ and let $\omega = (\omega_e)_{e \in E_d}$ be a random environment governed by a family of i.i.d. non-negative random variables defined on a probability space $( \Omega, \PP)$, which we call \textbf{random conductances} and also represent by $\omega_{xy}$. The expectation with respect to $\PP$ is denoted by $\E$. We assume that $\PP (\omega_e > 0 ) > p_c (d)$ where $p_c (d)$ is the critical parameter for Bernoulli bond percolation on $\Z^d$. Then, by virtue of standard percolation theory (cf. \cite{G}), there exists a unique infinite connected cluster $\mathscr{C}$ along positive conductances and on the other hand the set $\Z^d \setminus \mathscr{C}$ is a union of finite clusters that we commonly call \textit{holes}. 

Now, for any realization of the environment $\omega \in \Omega$, set
\begin{equation}
\label{tp}
\pi ( x ) 
= \sum_{y:| x - y | =1} \omega_{xy}, \qquad P_\omega (x, y) = \frac{\omega_{xy}}{\pi (x)} 
\end{equation}
Then, let $X = (X_n)$ be the discrete-time nearest-neighbor random walk in the random environment $\omega$ with jumping probabilities $P_\omega (x, y)$.
Let $P^x_\omega$ denote the \textit{quenched} probability law for $X$, started from $x$, and $E^x_\omega$ the associated expectation. This is obviously a reversible Markov chain with respect to $\pi$.

The \textit{equivalent} continuous-time version of $X$ waits on each site an exponential time of constant parameter $1$; the so-called constant speed random walk (CSRW) that we also denote (abusively) by $X$ with a time index $t$. $X$ is associated with the generator $P_\omega - I \eqqcolon \mathcal L$ defined on the Hilbert space $L^2 (\pi)$ which is equipped with the scalar product given by
$$
\langle f, g \rangle = \sum_{x \in \Z^d} f (x) g (x) \pi (x)
$$
Furthermore, the Dirichlet form associated with the continuous-time $X$ is defined by
$$
\mathcal E (f, f) = \langle - \mathcal L f, g \rangle = \frac{1}{2}\, \sum_{\{x, y\}} (f (x) - f (y))^2 \omega_{x y}
$$ 
 
But, if we let the random walk waits on each site $x$ an exponential time of parameter $\pi (x)$, we get the variable speed random walk (VSRW) which is denoted by $Y$; it is associated with the discrete Laplace generator $\pi \mathcal L$ and is reversible with respect to $1$.
In fact, $X$ and $Y$ are time-changes of each other. Indeed, define the additive functional 
$$
A (t) =  \int^t_0 \pi ( Y_u )\, {\rm d} u
$$
and let $A^{- 1} (t) = \inf \{s \ge 0 :\, A (s) > t\}$, that is, its right-continuous inverse.
Then, the stochastic process defined by $Y_{A^{-1} (t)}$ has the same law as $X$. Hence, we can consider that $X_t = Y_{A^{- 1} (t)}$ and with a slight abuse of notation, we also use $P^x_\omega$ to represent $Y$'s distribution.

The heat kernel associated with $X$ is defined by the quantity :
$$
P^t_\omega (x, y) = \frac{P^x_\omega (X_t = y)}{\pi (y)}
$$
For $Y$, this is simply $P^x_\omega (Y_t = y)$.
Besides, $P^o_\omega$ averaged on the environment is called the \textit{annealed} law which is given by
$$
\overline\PP (\,\cdot\,) = \int P^o_\omega (\,\cdot\,) \, {\rm d} \PP ( \omega )
$$

The RCM has witnessed an intense research activity in recent years where the main questions were to know the asymptotic behavior of the transition probabilities and the validity of the CLT in both cases, quenched and annealed  (see \cite{M-review} or \cite{K} for surveys).
Under the only condition of percolating positive conductances, the quenched CLT has been shown to be true (cf. \cite{ABDH}). However, things can become locally complicated because of ``traps", which may be due to either very small or very large conductances. This is reflected on the heat kernel behaviors, which has been studied in a large number of articles such as \twocite{FM}{BBHK}, \twocite{B1}{B2}, \twocite{UCLA-team}{BiBo} and \twocite{BKM}{B-SPL}.
In fact, we say that $X$ or $Y$ have a \textbf{normal} behavior if they behave as a symmetric simple random walk, in particular, if the return probability at time $t$ is of order $t^{- d/2}$, otherwise, the behavior is said to be \textbf{anomalous}. 

Henceforth, $c$ denotes a generic constant and we adopt for simplicity the following \textit{deterministic} notations : 
$f  \lesssim g$
to mean that $f = O (g)$. We write $f \simeq g$ when both $f \lesssim g$ and $f \gtrsim g$ hold. Plus, we use the usual notation $f \sim g$ to mean that $f/g = 1 + o (1)$. 

\section{\bf Results}

We address the question of slowing down of random walks in random environments of conductances and our main results concern the model with i.i.d. conductances in $[0, 1]$ and a law verifying the following polynomial condition with a parameter $\alpha > 0$,
\begin{equation}
\label{LP}
\tag{LP}
\PP ( \omega_e \le u ) \sim u^{- \alpha}, \quad u \to 0
\end{equation}

In a series of articles \treecite{B1}{B2}{BKM} and \cite{B-SPL}, we showed that the transition from normal to anomalous behavior occurs in this important model.
Indeed, under \eqref{LP} and using an original trapping method that has been confirmed in \cite{UCLA-team}, we obtained interesting anomalous return probability estimates in \cite{B1}, which were slightly improved in \cite{B-SPL} in the following form,
\begin{equation}
\label{SPL}
P^{2 n}_\omega (o, o) \gtrsim n^{- (2 + \varepsilon)}\, \Theta_n^{d - 1}
\end{equation}
where $\Theta_n$ is the probability of a special configuration of the environment that we call \textbf{trap} of depth $n$, that is,
\emph{a central \underline{trapping edge} $\mathbf{e} \in E_d$ with $\omega_{\mathbf{e}} \ge \theta$ for some positive $\theta < 1$, surrounded  by $4d - 2$ edges with $\omega_{e} \in [1/n,2/n]$.} 
Note that by \eqref{LP}, we have 
$$
\Theta_n \simeq n^{- (4 d - 2) \alpha}
$$

On the other hand, we have (see. \cite{BKM}) a normal heat kernel on all the sites of a diffusive-scaled box, i.e. $\sup_{x, y \in B_{\sqrt{t}}} P^t_\omega (x, y) \lesssim t^{- d/2}$ when $\alpha > \frac14\, d/(2d - 1)$, which yield a local-CLT for $X$, and also for $Y$ when $\alpha > 1/4$.    

Besides, there are general upper bounds that we can find in \cite{BBHK} : for percolating i.i.d. conductances in $[0, 1]$, for a.e. $\omega$ with $o \in \mathscr C$, we have 
\begin{equation*}
\label{trans}
P^{2 n}_\omega (o, o) \lesssim
\begin{cases}
n^{-d/2},\qquad & d = 2, 3
\\
n^{-2}\log n,\qquad & d=4
\\
n^{-2},\qquad & d \ge 5
\end{cases}
\end{equation*}
In the same time, the Markov property with the reversibility of $X$ and $\pi (x) \le 2 d$  combined with Cauchy-Schwarz inequality and the CLT $($see \cite{ABDH}$)$ give 
\begin{equation}
\label{gnlb}
P^{o}_\omega (X_{2 n} = o)
\ge 
\frac{\pi(o)}{2d}\, \frac{P^ o_\omega (|X_n| \le \sqrt{n})^2}{|\{\vert x\vert\leq \sqrt{n}\}|}
\gtrsim
\pi (o)\, n^{- d/2}
\end{equation}
Then, for $d = 2, 3$, the return probability is still normal. But, in $d = 4$, we proved in \cite{BiBo} that the extra-logarithmic term in the general upper bound expresses a real phenomenon and here is the lower bound,
\begin{equation}
\label{bibo}
P^{o}_\omega (X_{2 n} = o) \geq C (\omega)\, \frac{\log n}{n^2} \,\Theta_n^2
\end{equation}
for a constructed random environment of i.i.d. conductances in $[0, 1]$ such that  $\Theta_n^2 \ge \lambda_n^{- 1}$ where $\lambda_n$ grows to infinity as slowly as we want along a deterministic integers sequence. However, we have conjectured for $d = 4$ that in a general random environment of i.i.d. conductances in $[0, 1]$,
\begin{equation}
\label{C12}
\frac{n^2}{\log n}\, P^{2 n}_\omega (o, o) \xrightarrow[n \to \infty]{} 0, \qquad \PP-a.s.
\end{equation}
which has in the meantime been proved in \cite{UCLA-team}.

Our main results show under the assumption \eqref{LP} that for $\alpha$ small enough, the correct decay order of the anomalous return probability can be read as :
\begin{center}
\textit{the probability of a trap $\times$ the probability to get in and out of it.}
\end{center}
Formally, we have :

\begin{theorem}
\label{th1}
In a random environment on $\Z^d$ of i.i.d. random conductances in $[0, 1]$ satisfying \eqref{LP}, for $d \ge 2$, for any $\varepsilon > 0$, we have on the one hand, if $\alpha < 1/[2 (4 d - 2)]$,   
\begin{equation}
\label{A+}
P^{2 n}_\omega (o, o) \gtrsim n^{- (2 + \varepsilon)} \, \Theta_n, \quad \PP-a.s.,
\end{equation} 
and on the other hand, we have 
\begin{equation}
\label{A++}
\PP \left(P^{2 n}_\omega (o, o) \lesssim n^{- 2 + \varepsilon} \, \Theta_n \vee n^{- d/2}\right)
\xrightarrow[n \to \infty]{} 1
\end{equation}
\end{theorem}

\begin{remarks}
The bounds that appear in \eqref{A++}, \eqref{A+} and \eqref{gnlb} first show a normal behavior for $d = 4$, which answers the question $(2)$ in \cite[Remark~1.3]{BiBo} and confirms \eqref{C12}. But, for $d \ge 5$, the lower bound is anomalous if $\alpha$ is less than
\begin{equation*}
\alpha_c = \frac14\, \frac{d - 4}{2 d - 1}
\end{equation*}
Furthermore, the upper bound is normal if $\alpha > \alpha_c$.
Then, the right question is $:$ what is the critical value for $\alpha$, is it $\alpha_c$ ?!

On the other hand, $Y$ can be trapped in a different simpler configuration $:$ a given point $x$, \textbf{the trapping site}, surrounded by conductances of order $t^{-1}$; the probability of such a trap is $\Theta^*_t \simeq t^{- 2 d \alpha}$. 
If the random walk is held in this configuration for a time larger than $t$, we get 
$$
P^x_\omega (Y_t = x) 
\ge e^{- 2 d}
$$
Then, the same strategy adopted here provides for $\alpha < 1/(4 d)$ a lower bound of the form 
\begin{equation}
\label{A+Y}
P^{t}_\omega (o, o) \gtrsim t^{- (2 + \varepsilon)} \, \Theta_t^*
\end{equation}
which is anomalous when $\alpha$ is smaller than
\begin{equation*}
\alpha_v = \frac{1}{4}\, \frac{d - 4}{d}
\end{equation*}

Conversely, the same method developed for \eqref{A++} with a necessary change at Lemma~\ref{CeN} $($see Remark~\ref{CeNY}$)$ gives us an upper bound such as
\begin{equation}
\label{A++Y}
\PP \left(P^{t}_\omega (o, o) \lesssim n^{- 2 + \varepsilon} \, \Theta_t^* \vee t^{- d/2}\right)
\xrightarrow[n \to \infty]{} 1
\end{equation}
which is normal if $\alpha > \alpha_v$.
Therefore, is $\alpha_v$ the critical value for $Y$ ?!

In a second opposite situation, one can take conductances in $[1, \infty)$ satisfying a law with polynomial tail at infinity, i.e.
\begin{equation}
\label{UP}
\tag{UP}
\PP ( \omega_e > u ) \sim u^{- \alpha}, \quad u \longrightarrow \infty
\end{equation}
In this case, $Y$ is always normal $($cf. \cite{ABDH}$)$ but $X$ has an anomalous convergence when $\alpha < 1$ $($cf. \cite{BC}$)$. Actually, the same phenomenon occurs for $X$ because of traps which for this type of conductance are defined as follows $:$ \textit{a trap of depth $t$ is constituted of a central trapping edge $\mathbf e$ with $\omega_{\mathbf e} >  t$, surrounded by $4 d - 2$ incident edges with conductances in $[1, \theta]$ for some fixed $\theta > 1$}.
By using  \eqref{UP}, the probability of such a configuration is $\Theta_t \simeq t^{- \alpha}$.

Explicitly, let $X$ be the CSRW in a random environment on $\Z^d$ of i.i.d. conductances in $[1, \infty)$ obeying \eqref{UP}. 
Then, for $d \ge 3$, for any $\varepsilon > 0$, for a.e. $\omega$, we have 
\begin{equation}
\label{arp}
P^o_\omega (X_t = o) 
\begin{cases}
= o\, (t^{- 1}) & \forall \alpha > 0
\\
\simeq t^{- d/2} & \alpha > 1 
\\
\gtrsim
\pi (o)\, t^{- (1 + \varepsilon)}\, \Theta_t & \alpha <  1/2 
\end{cases}
\end{equation}

For $\alpha \in [1/2, 1]$, the question remains open but we believe in view of \cite{C} that we still have an anomaly. Concerning $d = 2$, we also expect by \cite{C} an anomaly that would be expressed with an extra-logarthmic term. 
In the same time, the present proof that gives \eqref{A++} unfortunately does not apply to 
\eqref{UP}.  
\end{remarks}

The rest of the paper consists of tree sections that give proofs of \twoeqref{A+}{A++} and \eqref{arp}.



\section{\bf The lower bound}\label{}

Henceforth, we assume that the conductances are i.i.d. random variables in $[0, 1]$ and satisfy \eqref{LP}. Thus, the conductances are positive a.s. and $\scrC \stackrel{\rm a.s.}{=} \Z^d$. 

Our proof of \eqref{A+} basically optimizes the trapping strategy used in \twocite{B1}{B-SPL} which consists in showing that the random walk spends most of its time in traps as defined above.
To this end, set $B_k = [- k, k]^d \cap \Z^d$ and let $\partial B_k = B_k \setminus B_{k - 1}$, that is, the inner boundary of $B_k$. Call $H_k$ the hitting time of $\partial B_k$, i.e. $H_k = \inf \{t \ge 0 : \, X_n \in \partial B_k\}$.
Then, define $V_k$ as the event that at $H_k$, the random walk meets a trap of depth $n$ which is located outside of $B_k$ and call $\textbf{e}_k = \{\textbf{v}_1, \textbf{v}_2\}$ its trapping edge (see \twocite{B1}{B-SPL} for more details).

Let $m = \lceil \Theta_n^{- 1} \, n^{\varepsilon} \rceil$ with an arbitrary $\varepsilon  > 0$ and note the following first key step that $X$ meet a.s. a trap before exiting $B_{m}$ :

\begin{lemma}
\label{trap-l}
For a.e. $\omega$,  for all $n$ large enough, we have
\begin{equation}
\label{Lm}
P^{o}_\omega \left(\sum^m_{k = 1} \1_{V_k} = 0\right) \le e^{- \frac15\, n^{\varepsilon}}
\end{equation}
\end{lemma}

\begin{proof}
Following a similar argument for the proof of  \cite[Lemma~2.3]{B-SPL}, let $\Lambda_{m}$ be the event that there is a trap at one of the locations $X_{H_{k}}$, i.e.
\begin{equation}
\label{T_N}
\Lambda_{m} = \bigcup_{k= 1}^{m} V_{k}
\end{equation}  
and note that $H_k < \infty$ a.s., since $X$ is an irreducible Markov chain. Second, remark that the $V_{k}$ are by construction $\PP \times P^ o_\omega$-independent if $k \in 3 \N$. Plus, we have by invariance of the environment $\overline\PP \left(V_k\right) = \Theta_n $. Hence, we obtain that
$$
\overline\PP \left(\Lambda_{m}^c\right) \le \big(1 - \Theta_n\big)^{[m/3]}
$$

Now, use Markov inequality to get
\begin{equation*}
\PP \big(P^ o_\omega (\Lambda_{m}^c\big) > e^{- \frac16\, m \Theta_n})
\le 
e^{\frac16\, m \Theta_n}\, \big(1 - \Theta_n\big)^{[m/3]}
\le 
e^{- \frac14\, m \Theta_n}
\lesssim e^{- \frac14\, n^{\varepsilon}}
\end{equation*}
Then, we deduce by Borel-Cantelli lemma that $\PP-$a.s.,
$$
P^ o_\omega (\Lambda_{m}^c) 
\le e^{- \frac15\, n^{\varepsilon}} 
$$
The claim follows.
\end{proof}

Nevertheless, the random walk remains confined in a diffusive-scaled box. Explicitly, by virtue of Carne-Varopoulos inequality (cf. \cite{MR}), the exit time $\tau_m$ from $B_m$ satisfies
\begin{equation}
\label{taum}
P^o_\omega (\tau_m \le n) \le 2 n\, m^{d - 1} e^{- m^2/n} + e^{- c n}
\end{equation}
Thus, we need $\Lambda_m$ from \eqref{T_N} to occur for $m < \sqrt{n}$ and the following hitting time estimate is necessary.

\begin{lemma}
\label{H_N}
We have $\PP-$a.s.,
\begin{equation}
E^ o_\omega (H_m) \lesssim m^{2}\, (\log m)^{2 d}
\end{equation}
\end{lemma}

The proof of this expectation requires a time change on $X$.
First, pick a $\theta > 0$ such that $\PP (\omega_e > \theta ) > p_c (d)$. Standard percolation guarantees (cf. \cite{G}) that there exists almost surely a unique infinite cluster of bonds with conductances larger than $\theta$ that we denote by $\scrC^\theta$, the so-called \textit{strong component} which also represents the set of all vertices in its bonds. On the other hand, we have the cluster $\mathscr H^\theta \coloneqq \Z^d\setminus \mathscr C^\theta$ which is a union of finite \textit{holes}. Note that if $\mathscr H_o^\theta$ is the weak component incident to $o$ or contains $o$, we have from \cite[Lemma~3.1]{BP} that for $\theta$ small enough, there exists a constant $c > 0$ such that 
\begin{equation}
\label{diamH}
\PP ( \diam\, {\scrH}_o^\theta > m ) \lesssim e^{- c m}
\end{equation}
Here ``$\diam$" is the diameter in the $|\,\cdot\,|_\infty-$distance on $\Z^d$. 

Then, let $X^\theta = (X^\theta_l)$, the \textit{coarse grained} random walk which records the successive visits of $X$ to $\mathscr C^\theta$. Explicitly, set $T_0 = 0$ and let $(T_l)_{l \ge 1}$ be the successive time lapses that $X$ spends in holes before returning to $\mathscr C^\theta$, i.e. 
$$
T_l = \inf \{n > 0 : X_{T_l + n} \in \mathscr C^\theta\}
$$
Write $S_l = \sum^l_{i = 1} T_l$ and define
\begin{equation}
\label{Xthl}
X^\theta_l = X_{S_l}
\end{equation}
Thus, $X^\theta$ is a Markov chain with transition probabilities given by $P^x_\omega (X_{T_1} = y)$ for $x, y \in \mathscr C^\theta$ with an invariant and reversible measure $\pi$ restricted to $\mathscr C^\theta$; the conductances are induced by the transition probabilities and can join sites separated by a distance larger than $1$ through a hole.

The random walk $X^\theta$ has a normal behavior : we have (see for eg. \cite[Lemma~3.2]{BBHK}) that
for a.e. $\omega$, for $x, y \in \scrC^\theta$, for $l$ larger than a random value $L (\omega)$ that depends on the isoperimetry of $\mathscr C^\theta$,
\begin{equation}
\label{CGXN}
P^x_\omega (X^\theta_l = y) \lesssim \pi (x)\, l^{- d/2}
\end{equation}

\begin{proof}[Proof of Lemma~\ref{H_N}]
Let $\tau_m^\theta$ be the exit time of $X^\theta$ from $B_m$. Then, as a consequence of \eqref{CGXN}, we have from \cite[Lemma~2.1]{B-SPL} that
\begin{equation}
\label{tthm}
E^o_\omega (\tau^\theta_m) \lesssim m^2
\end{equation}

Now, observe that
$$
H_m \le \sum_{l = 0}^{\tau_m^\theta} T_l
$$
which implies by the Markov property, 
\begin{equation}
\label{Hm0}
E^ o_\omega (H_m) 
\le 
E^ o_\omega \bigg( \sum_{l = 0}^{\tau_m^\theta}  E^{X^\theta_l}_\omega (T_1)\bigg)
\le 
E^ o_\omega \big(\tau_m^\theta\big)\, \sup_{x \in \scrC^\theta_m} E^{x}_\omega (T_1)
\end{equation}
where $\scrC^\theta_m \coloneqq \scrC^\theta \cap B_m$.

At the same time, \cite[Lemma~3.8]{BBHK} tells us that for $\theta$ small enough, for a.e. $\omega$ and for all $x \in \scrC^\theta$,
\begin{equation}
\label{Tl+}
E^x_\omega (T_1) \le \frac{4 d}{\theta}\, |\mathscr H_x^\theta|
\end{equation}
Besides, \eqref{diamH} and Borel-Cantelli lemma yield that 
$$
\sup_{x \in \scrC^\theta_m} |\mathscr H_x^\theta| \le (\log m)^{2 d} 
$$ 
This and \eqref{tthm} give the desired estimate by substituting in \eqref{Hm0}.
\end{proof}

Finally, we claim the following :

\begin{lemma}
\label{XinS}
If $\alpha (4 d - 2) + \varepsilon < 1/2$, we have for a.e. $\omega$, for $k = 1, \ldots, m$,  
\begin{equation*}
\label{stayintrap}
P^ o_\omega \big(X_n \in {\rm  \mathbf{e}}_k \mid V_{k}\big) \gtrsim  n^{-1}
\end{equation*}
\end{lemma}

\medskip
\begin{proof}
Let $n$ be large enough such that Lemma~\ref{trap-l} holds. Suppose $V_k$ occurs, then define the event 
$$
W_k = \bigcap^{H_{k} + n}_{j = H_{k} + 1} \big\{ X_j \in \mathbf{e}_k\big\}
$$
This is the event that $X$ spends a time $n$ on the trapping edge $\mathbf{e}_k$. Moreover, note that the probability to cross $\mathbf{e}_k$ is larger than $1/(2 d n)$. Hence, we obtain by the Markov property
\begin{equation}
\label{Dk}
P^{X_{H_k}}_\omega \big(W_k \mid V_k\big) 
\ge \frac{1}{2 d n}\, \left( 1 - \frac{2 d - 1}{\theta n}\right)^{n} 
\gtrsim 
n^{- 1}
\end{equation}

Now, remark that on $V_k$,
\begin{equation*}
\{H_{k}< n\} \cap  W_k \subset \{X_n \in \mathbf{e}_k\}
\end{equation*} 
which, by the Markov property and \eqref{Dk}, implies that
\begin{eqnarray*}
P^ o_\omega \big(X_n \in \mathbf{e}_k \mid V_k\big) 
\ge
E^ o_\omega \left(\1_{\{H_k < n\}}\, P^{X_{H_k}}_\omega (W_k \mid V_k)\right) 
\gtrsim
n^{-1}\, P^ o_\omega \big(H_{k} < n\big)
\end{eqnarray*}
But, by Lemma~\ref{H_N} and Markov inequality with the condition $\alpha (4 d - 2) + \varepsilon < 1/2$, we have
\begin{equation*}
P^ o_\omega \left(H_{k} \ge n\right)  \le P^ o_\omega \left(H_{m} \ge n\right) 
\lesssim m^{2} (\log m)^{2 d} n^{- 1}
\xrightarrow[n \to \infty]{} 0
\end{equation*}
\end{proof}

Therefore, we are ready to finalize the demonstration of the anomalous lower bound.

\begin{proof}[Proof of \eqref{A+}]
First, use that $X$ cannot be in different $\mathbf{e}_k$ in the same time and the Markov property to obtain that 
\begin{align*}
\label{2Vk}
P^ o_\omega (X_{2 n} = o)
&\ge
\sum^{m}_{k = 1} P^o_\omega (X_{n} \in \mathbf{e}_k, X_{n} =o) 
= \sum^{m}_{k = 1} \sum_{i = 1, 2} E^o_\omega \left(\1_{\{X_{n} = \textbf{v}_i\}} P^{\textbf{v}_i}_\omega (X_{n} = o)\right)
\end{align*}
which, by the reversibility and $\pi (\textbf{v}_i) \le 2 d$, plus Cauchy-Schwarz inequality, is larger than
\begin{align*} 
\frac{\pi (o)}{2 d}\,\sum^{m}_{k = 1} \sum_{i = 1, 2} P^o_\omega (X_ {n}= \textbf{v}_i)^2
\ge 
\frac{\pi (o)}{4 d}\,\sum^{m}_{k = 1} P^o_\omega (X_{n} \in \mathbf{e}_k)^2
\ge 
\frac{\pi (o)}{4 d m}\, \left(\sum^{m}_{k = 1} P^o_\omega (X_{n} \in \mathbf{e}_k)\right)^2
\end{align*}
Next, take this sum squared and observe that
\begin{eqnarray*}
\sum^{m}_{k = 1} P^o_\omega (X_{n} \in \mathbf{e}_k)
&=&
\sum^{m}_{k = 1} P^{o}_\omega (X_{n} \in \mathbf{e}_k \mid V_k)\, P^{o}_\omega (V_k)
\\
&\gtrsim&
n^{- 1}\, E^{o}_\omega \left(\sum^m_{k = 1} \1_{V_k}\right)
\\
&\gtrsim&
n^{- 1}\, P^{o}_\omega \left(\sum^m_{k = 1} \1_{V_k} \ge 1\right)
\end{eqnarray*}
Then, use Lemma~\ref{trap-l} and replace in the second development above to get the desired lower bound.
\end{proof}

\section{\bf The upper bound}\label{}

This section is devoted to demonstrating \eqref{A++} in two important parts which subtly combines the argument used in \cite{BBHK} with spectral analysis techniques that have already been employed in \twocite{FM}{BKM}. 
Of course, the assumption \eqref{LP} is from now on assumed.

Let $A_\theta (n)$ be the discrete-time additive functional that records the time spent in $\mathscr C^\theta$, i.e. $A_\theta (n) = \sum^n_{i = 0} \1_{\mathscr C^\theta} (X_i)$.
Furthermore, let $\xi > 0$ and set $N = n^{(1 + \xi)/2}$ for Carne-Varopoulos inequality reason that appears later. 

First, recall from \cite[Lemma~4.2]{BKM} that for any $\eta \in (0, 1)$, for $\theta$ small enough, we have for $\mathscr C^\theta_N = \mathscr C^\theta \cap B_N$,
\begin{equation}
\label{eta}
\PP (|\mathscr C^\theta_N| \le \eta |B_N|) \le e^{- c N}
\end{equation}
So it comes that
\begin{equation}
\label{eta+}
\mathbb E \left(\frac{|\overline{\mathscr C^\theta_N}|}{|B_N|}\right) \le 1 - \eta +  e^{- c N} 
\end{equation} 

Now, let $\varepsilon > 0$ and set $m = \Theta_n^{- 1} n^{2 - \varepsilon}$. Then, by virtue of the invariance of the environment law $\PP$ and by using that the return probability is decreasing (cf.  \cite[Lemma~3.9]{BBHK}) with Markov inequality, we obtain that
\begin{eqnarray}
\PP \left(P^o_\omega (X_{2 n} = o) > m^{- 1}\right) 
&=& 
|B_N|^{- 1}\, \sum_{x \in B_N} \PP \left(P^x_\omega (X_{2 n} = x) > m^{- 1}\right) 
\nonumber
\\
&\le& 
|B_N|^{- 1} \mathbb E \bigg(\sum_{x \in \mathscr C^\theta_N} P^x_\omega (X_{2 n} = x) > m^{- 1}\bigg) 
 + \mathbb E \bigg(\frac{|\overline{\mathscr C^\theta_N}|}{|B_N|}\bigg)
\nonumber
\\
\label{S1}
&\stackrel{\eqref{eta+}}{\le}& 
|B_N|^{- 1}\, \frac{m}{n}\, \mathbb E \bigg(\sum_{x \in \mathscr C^\theta_N} \sum_{m = n + 1}^{2 n} P^x_\omega (X_{m} = x)\bigg)
\\
\nonumber
&& \qquad \qquad + 1 - \eta +  e^{- c N}
\end{eqnarray}

Hence, we have to deal with the double sum in \eqref{S1}. 
Therefore, recall from \eqref{Xthl} the coars-grained random walk $X^\theta$ and observe that
\begin{eqnarray}
\sum_{m = n + 1}^{2 n} P^x_\omega (X_{m} = x)
&=&
\sum_{m = n + 1}^{2 n} \sum^{m}_{l = 1} P^x_\omega (X^\theta_{l} = x, S_{l} = m)
\nonumber
\\
&=& 
\sum^{n + 1}_{l = 1} \sum^{2 n}_{m = n + 1}  P^x_\omega (X^\theta_{l} = x, S_{l} = m)
\nonumber
\\
&&
\qquad + \sum^{2 n}_{l = n + 2} \sum^{2 n}_{m = l}  P^x_\omega (X^\theta_{l} = x, S_{l} = m)
\nonumber
\label{enl0}
\\
&\le&
\sum^{[\varepsilon n]}_{l = 1} P^x_\omega (X^\theta_{l} = x,  n \le S_{l} \le 2 n)
\\
\label{enl}
&&
\qquad +
 \sum^{2 n}_{l =[\varepsilon n]} P^x_\omega (X^\theta_{l} = x, S_{l} \ge n)
\end{eqnarray}

On the one hand, the sum in \eqref{enl0} is less than
\begin{eqnarray}
\label{n*}
\sum^{[\varepsilon n]}_{l = 1} P^x_\omega (X^\theta_{l} = x, A_\theta (2 n) = l)
\le  
P^x_\omega (X_{2 n} = x, A_\theta (2 n) \le \varepsilon n) 
\end{eqnarray}

On the other hand, consider the terms of the sum in \eqref{enl} and  recall \eqref{CGXN}. As a consequence, \mbox{we have the following :}
$\PP-$a.s., for $x \in \mathscr C^\theta$, for every $n \ge 1$, we have for $l$ large enough, 
\begin{equation}
\label{XSl}
P^x_\omega (X^\theta_l = x, S_l \ge n) \lesssim \frac{l^{1 - d/2}}{n}
\end{equation}

This is similar to \cite[Proposition~3.5]{BBHK} and here we quickly reproduce the proof. First remark that by reversibility, we have
\begin{eqnarray*}
P^x_\omega (X^\theta_l = x, S_l \ge n)
\le
2\, P^x_\omega (X^\theta_l = x, S_{[l/2]} \ge n/2) 
\end{eqnarray*}
Take this doubled probability, use Markov inequality and condition on the position of $X^\theta$ at the times before and after $T_i$ to get
$$
P^x_\omega (X^\theta_l = x, S_{[l/2]} \ge n) 
\le \frac{2}{n}\, \sum^{[l/2]}_{i = 1} \sum_{v, w} P^x_\omega (X^\theta_{i - 1} = v) E^v_\omega (T_1; X^\theta_i = w) 
P^w_\omega (X^\theta_{l - i} = x)
$$
But $X^\theta$ is normal (see \eqref{CGXN}), which, by reversibility and $l - i \ge l/2$, yields that 
$$
P^w_\omega (X^\theta_{l - i} = x) \lesssim l^{- d/2}
$$
Thus, it remains 
$$
\sum_w E^v_\omega (T_1; X^\theta_i = w) = E^v_\omega (T_1) \stackrel{\eqref{Tl+}}{\le} \frac{4 d}{\theta}\, |\mathscr H_v^\theta| 
$$
where $\mathscr H_v^\theta = \mathscr H^\theta (v)$ is the weak component incident to $v$. 
Then, it comes that
$$
P^x_\omega (X^\theta_l = x, S_{[l/2]} \ge n) \lesssim \frac{l^{- d/2}}{n}\, E^x_\omega \left(\sum^{[l/2] - 1}_{i = 0} |\mathscr H^\theta ({X^\theta_i})|\right) 
$$

In the same time, we have by ergodicity (see  \cite[Lemma~3.7]{BBHK}) that 
$$
l^{- 1}\, E^x_\omega \left(\sum^{[l/2] - 1}_{i = 0} |\mathscr H^\theta ({X^\theta_i})|\right) \xrightarrow[l \to \infty]{\rm a.s. \,\&\, L^1} \mathbb E_\theta ( |\mathscr H_o^\theta|)
$$
where $\mathbb E_\theta$ is the expectation with respect to the measure
$$
\frac{\pi (o)}{\mathbb E (\pi (o) \mid o \in \mathscr C^\theta)}\, \mathbb P ({\rm d}\, \omega \mid o \in \mathscr C^\theta)
$$
The limit is finite by \eqref{diamH}. Thus, the desired bound \eqref{XSl} follows.

Therefore, \eqref{n*} and \eqref{XSl} give for $d \ge 4$,
\begin{equation}
\label{n0}
\sum_{x \in \mathscr C^\theta_N} \eqref{enl0} + \sum_{x \in \mathscr C^\theta_N} \eqref{enl} 
\lesssim
\sum_{x \in \mathscr C^\theta_N} P^x_\omega (X_{2 n} = x, A_\theta (2 n) \le \varepsilon n) 
+ |B_N|\, n^{1 - d/2}
\end{equation}
However,  by Carne-Varopoulos inequality (see \eqref{taum}), the exit time $\tau_N$ from $B_N$ with $N = n^{(1 + \xi)/2}$ satisfies,
\begin{equation}
\label{tauN}
P^x_\omega (\tau_N \le n)) \le e^{- c n^\xi}
\end{equation}
This means that our focus should be on the sum 
\begin{equation*}
\sum_{x \in \mathscr C^\theta_N}P^x_\omega (X_{2 n} = x, A_\theta (2 n) \le \varepsilon n, \tau_N > n) \eqqcolon K_N
\end{equation*}
In other words, we consider $X$ with Dirichlet conditions outside $B_{2 N}$, that is, we follow the progress of the random walk until it leaves $B_{2 N}$ for the first time where it is killed. 

Here we move on to the second part of our demonstration, which will be essentially based on spectral analysis method. 
For reasons of technical simplicity, from now on, we switch from discrete-time to continuous-time by replacing $n$ by $t$. This is indeed possible because of the following. If $N (t)$ is the Poisson process that counts the number of jumps of $X$, up to time $t$, then we have for any $x$ and all $t$, 
\begin{equation}
\label{Nt}
P^x_\omega (N (t) \in [t/2, 2 t]) \lesssim (2/e)^{t/2}
\end{equation}

Recall that $\mathcal L$ is the generator of $X$ and let $L_{N}^2 (\pi)$ be the restricted space of functions defined on $B_{2 N}$ and null elsewhere.
Then, call $\mathcal L_{N}$ the restriction of $\mathcal L$ on $L_{N}^2 (\pi)$ and remark that $- \mathcal L_{N}$ is a positive symmetric operator.

The idea is to control $K_N$ by dividing the environment into ``good'' an ``bad'' edges. Explicitly, for bond percolation with parameter $p = 1 - q = \PP (\omega_e > u)$, write $\mathscr C_u (e)$ for the connected cluster that contains $e$ and note the following technical percolation lemma.  

\begin{lemma}
\label{CeN}
For $u$ small enough, we have
\begin{equation}
\label{CeN+}
\PP (|\mathscr C_u (e)| < \infty) \lesssim u^{(4 d - 2) \alpha}
\end{equation} 
\end{lemma}

\begin{remark}
\label{CeNY}
For $Y$, we would need to cut a site instead of an edge. Thus, if $\mathscr C_{u} (x)$ is the connected cluster that contains $x$,  \cite[Lemma~4.3]{BKM} gives
$$
\PP (|\mathscr C_u (e)| < \infty) 
\lesssim u^{2d\, \alpha}
$$
\end{remark}

\begin{proof}
By \cite[Lemma~4.3]{BKM} and \eqref{LP}, we have for $q$ small enough,
$$
\PP (|\mathscr C_u (e)| < \infty) \lesssim q^{(4 d - 2)} \lesssim u^{(4 d - 2) \alpha}
$$ 
\end{proof}

Now, choose $u = t^{\delta - 1}$ with $\delta \in (0, 1)$ and define $\Omega_N = \{\forall e \in B_N :\, |\mathscr C_{u} (e)| = \infty\}$. Remark that $\Omega_N \xrightarrow[]{\rm a.s.} \Z^d$ as $u \downarrow 0$. Then, consider first $K_N$ on $\Omega_N$ and observe that for any $\lambda > 0$ and by Markov inequality, we obtain
\begin{equation*}
K_N \le e^{\varepsilon \lambda t} 
\sum_{x \in \mathscr C^\theta_N} E^x_\omega \left(\delta_x (X_t)\, e^{- \lambda A_\theta (t)}; \tau_N > t\right)
\end{equation*}
where $A_\theta (t) = \int^t_0 \1_{\mathscr C^\theta} (X_s)\, {\rm d} s$. 
Hence, we are concerned with the semigroup $R_N = (R^t_N)$ defined by the Feynman-Kac formula :
$$
R^t_N f (x) = E^x_\omega \left(f (X_t) e^{- \lambda A_\theta (t)}; \tau_N > t\right)
$$

In $L^2_N (\pi)$, $R_N$ is generated by the self-adjoint operator $\mathcal G^\lambda_N = \mathcal L_N - \lambda \mathcal{M}$ where $\mathcal M$ operates so that $\mathcal M f (x) = \1_{\mathscr C^\theta} (x) f (x)$. 
Besides, $- \mathcal G_N^\lambda$ is a positive symmetric operator. Let us then call $\lambda_i, i = 1,\ldots, |B_{2 N}|$, the sequence of its eigenvalues labeled in increasing order which are obviously associated with a set of eigenvectors duly normalized in $L_N^2 (\pi)$. 

Therefore, we arrived on $\Omega_N$ at the inequality 
\begin{equation}
\label{SR}
K_N \le e^{\varepsilon \lambda t} \sum_{x \in \mathscr C^\theta_N} R^t_N \delta_x (x)  
\end{equation}
which leads to the following result.

\begin{proposition}
For any $\epsilon > 0$, for $t$ large enough, we have
\begin{equation}
\label{SNO}
\mathbb \E (K_N; \Omega_N) 
\le e^{- c t^\epsilon}
\end{equation}
\end{proposition}

\begin{proof}
Consider the sum in \eqref{SR} and observe that by using the trace formula of $R_N$,
\begin{equation*}
\sum_{x \in \mathscr C^\theta_{N}} R^t_N \delta_x (x) \le \sum_{x \in B_{N}} R^t_N \delta_x (x)
\le \sum_{i = 1}^{|B_{2 N}|} e^{- \lambda_i t}
\le |B_{2 N}| e^{- \lambda_1 t}
\end{equation*}
Thus, we clearly have to estimate the spectral gap $\lambda_1$ which is given by the following variational formula on $L_N^2 (\pi)$,
\begin{equation}
\label{l1}
\lambda_1 = \inf_{f \neq 0} \frac{\langle - \mathcal G_N f, g \rangle}{\langle f, g \rangle} 
\end{equation}

We need to impose two further conditions on $\Omega_N$. Let $m = [t^\epsilon]$ for an arbitrary $\epsilon > 0$ and set $B_m (z) = (2 m + 1) z + B_m$ for $z \in \Z^d$. Hence, $(B_m (z))_{z \in B_{2 N}}$ forms a partition of $B_{2 N}$. Call $D_N^1$ the event that $|\mathscr C^\theta \cap B_m (z)| > 2/3\, |B_m|$ for all $z \in B_{2 N}$ and let $D_N^2$ be the one that $|\mathscr H_z^\theta| \le N^\epsilon$ for all $z \in \mathscr C^\theta_{2 N}$.
Then, define
$$
\Omega_N^* = \Omega_N \cap D_N^1 \cap D_N^2
$$
We claim that
on $\Omega_N^*$ with $\lambda = t^{- 1 + \delta/2}$, we have for $t$ large enough,
\begin{equation}
\label{l1+}
\lambda_1 \ge t^{- 1 + \delta/2}
\end{equation}

To prove this, we are going to use path method. Remark first that the definition of $ \mathcal G_N$ and \eqref{l1} yield 
\begin{equation}
\label{sgvf}
\lambda_1 = \inf_{f \neq 0} \frac{\mathcal E (f, f) + \sum_{\mathscr C^\theta_{2 N}} f (x)^2 \pi (x)}{\sum_{x \in B_{2 N}} f (x)^2 \pi (x)}
\end{equation}

Next, let $x \in \mathscr H_{2 N}^\theta \coloneqq B_{2 N} \setminus \mathscr C^\theta_{2 N}$ and call $e = \{x, y\}$ the edge such that \mbox{$\omega_e = \max_{y \sim x} \omega_{x y}$}. Then, on $\Omega_N^*$, we clearly have a path $l (e, x^*)$ that connects $e$ to a site $x^*$ in $\mathscr C^\theta \cap B_m (z)$ along conductances larger than $t^{\delta - 1}$. We can choose this connection to be injective thanks to $D^1_N$.

Now, suppose that $l (e, x^*)$ starts at $y$, otherwise the calculation is easier. So, we have
\begin{equation*}
\label{dfe}
f (x) = f (x) - f (y) + \sum_{b\in l( e, x^*)} \text d f(b) + f(x^*),
\end{equation*}
which, by Cauchy-Schwarz inequality, implies
\begin{equation*}
\label{dfer}
f (x)^2 \leq 2 ( f (x) - f (y) )^2 + 4 \vert \ell(e, x^*)\vert \sum_{b \in l( e, x^*)} \text d f(b)^2+ 4 f (x^*)^2
\end{equation*}  
Multiply this inequality by $\pi (x)$ to obtain
\begin{equation}
\label{fx2}
f(x)^2\, \pi (x) \leq 4d ( f (x) - f (y) )^2\omega_e 
+ 8d |l (e, x^*)| \sum_{b \in l( e, x^*)} \text d f(b)^2+  \frac{8 d}\theta f (x^*)^2\,\pi (x^*)
\end{equation}
where we used that $\pi (x^*) \ge \theta$ and $\pi (x) \le 2d\, \omega_e \le 2d$.

But, $l( e, x^*)$ must be in a \textit{negligible} size hole $\mathscr H_{x^*}^\theta \subset \mathscr H_{2 N}$. Indeed, it comes from \eqref{diamH} that for any $\epsilon > 0$, we have $\PP-$a.s.,
\begin{equation}
\label{SNs}
|l (e, x^*)| \le |\scrH_{x^*}^\theta| \le N^\epsilon
\end{equation}

Then, on $\Omega_N^*$, \eqref{fx2} becomes 
$$
f (x)^2\, \pi (x) \leq 4d ( f (x) - f (y) )^2 \omega_e 
+ 8d N^\epsilon\, n^{1 - \delta} \sum_{b \in l ( e, x^*)} \text d f (b)^2 \omega_b +  \frac{8 d}\theta f (x^*)^2\,\pi (x^*)
$$
which, by summing over $x$, gives
$$
\sum_{x \in \scrH_{2 N}}  f (x)^2\, \pi (x) \leq 8d N^{2 \epsilon}\, t^{1 - \delta} \mathcal E (f, f) 
+ \frac{8 d}\theta \sum_{x\in \scrC_{2 N}^\theta} f (x)^2\,\pi (x)
$$
where we used that a bond can be repeated $N^\epsilon$ times and $x^*$ appears at most once. Hence, we get
$$
\sum_{x\in \scrC_{2 N}^\theta} f (x)^2\,\pi (x) \le \frac{16 d}\theta\, t^{1 - \delta + \epsilon (1 + \xi)}\mathcal E (f, f) + \frac{16 d}\theta  \sum_{x\in \scrC_{2 N}^\theta} f (x)^2\,\pi (x)
$$
Since $\epsilon$ and $\xi$ are arbitrary, this yields the claimed lower bound on the spectral gap by using \eqref{sgvf} with $\lambda = t^{- 1 + \delta/2}$.

Now, back to our primary objective, i.e. controlling $K_N$ on $ \Omega_N$, we get by \eqref{l1+} that
\begin{equation*}
\mathbb E (K_N; \Omega_N^*) \le |B_{2 N}|\, e^{\varepsilon \lambda t}\, \mathbb E (e^{- \lambda_1 t}; \Omega_N^*) 
\le |B_{2 N}| e^{- (1 - \varepsilon) t^{\delta/2}}
\end{equation*}
At the same time, we have by \eqref{eta} that for $\theta$ small enough,
\begin{equation*}
\PP \left(\overline{D_N^1}\right) \lesssim |B_{2 N}| e^{- c m}
\end{equation*} 
Besides, \eqref{diamH} gives
\begin{equation*}
\PP \left(\overline{D_N^2}\right) \lesssim |B_{2 N}| e^{- c N^{\epsilon/d}}
\end{equation*}
Then, we obtain
\begin{equation*}
\mathbb E (K_N; \Omega_N) 
\le 
|B_{2 N}| e^{- (1 - \varepsilon) t^{\delta/2}} + |B_{2 N}|^2 e^{- c t^{\epsilon/(2 d)}}
\end{equation*}
where we used at the end that $|K_N| \le |B_{2 N}|$. Since $\delta$ and $\epsilon$ are arbitrary, the desired inequality follows.
\end{proof}

Second, consider $K_N$ on $\Omega_N^c$ and let $\{\sigma_i, i = 1,\ldots, |B_{2 N}|\}$ be the set of eigenvalues of $- \mathcal L_N$ labeled in increasing order which are associated with the set of eigenvectors $\{\phi_i, i = 1, |B_{N}|\}$ with due normalization in $L_N^2 (\pi)$. 

\begin{proposition}
We have
\begin{equation}
\label{SNOc}
\mathbb E (K_N; \Omega_N^c) \lesssim |B_{2 N}|\,  t^{\varepsilon/2 - 1 - \alpha (4 d - 2)} 
\end{equation}
\end{proposition}

\begin{proof}
Remark that
$$
K_N \le \sum_{x \in \mathscr C^\theta_N} P^x_\omega (X_t = x, \tau_N > t)
$$
But, we have
$$
P^x_\omega (X_t = x, \tau_N > t) = \pi (x) \sum_{i = 1}^{|B_{2 N}|} e^{- \sigma_i t} \varphi_i (x)^2
$$
Besides, observe that
\begin{eqnarray*}
\sum_{x \in \mathscr C^\theta_N} \phi_i (x)^2
\le
\frac12 \sum_{\{x, y\} \in \mathscr C^\theta_N} (\phi (x) - \phi (y))^2
\le
\frac1\theta\, \mathcal E (\phi_i, \phi_i)
=
\frac{\sigma_i}\theta 
\end{eqnarray*}
Thus, it comes that
\begin{eqnarray*}
K_N
\le
2 d\, \sum_i e^{- \sigma_i t} \sum_x \phi (x)^2
\le 
\frac{2 d}\theta\, \sum_i \sigma_i e^{- \sigma_i t} 
\le \frac{2 d}\theta\, t^{- 1}\, \sum_i  e^{- \sigma_i t/2} 
\end{eqnarray*}
where we used at the end that $ \sigma_i e^{- \sigma_i t} \le t^{- 1} e^{- \sigma_i t/2}$. Hence, we obtain 
$$
\mathbb E (K_N; \Omega_N^c) \lesssim |B_{2 N}|\, t^{- 1} \PP (\Omega_N^c)
$$

It remains to estimate $\PP (\Omega_N^c)$. First, note that $|\mathscr C_{u} (e)| < \infty$ can only happen for $e \in \mathscr H^\theta_{2 N}$ since the edges in $\mathscr C^\theta_{2 N}$ percolate. So, let us partition $\Omega_N^c$ on the holes of $B_{2 N}$, that is, (with a slight notational abuse) for a hole $\mathscr H^\theta$ in $B_{2 N}$, the conditional probability that the edges in $\mathscr H^\theta$ are cut by conductances less than $t^{\delta - 1}$ is
\begin{equation*}
\PP (\Omega_N^c \mid \mathscr H^\theta) 
\stackrel{\eqref{SNs}}{\le} N^\epsilon t^{(\delta - 1) \alpha (4 d - 2)} \le t^{\varepsilon/2 - \alpha (4 d - 2)}
\end{equation*}
where we have chosen that $\epsilon (1 + \xi)/2 + \delta \alpha (4 d - 2) \le \varepsilon/2$. Therefore, we get
\begin{equation*}
\PP (\Omega_N^c)
= 
\mathbb E \bigg(\sum_{\mathscr H^\theta \subset \mathscr H_{2 N}^\theta} \PP (\Omega_N^c \mid \mathscr H^\theta) \1_{\mathscr H^\theta}\bigg)
\stackrel{\eqref{SNs}}{\le}
t^{\varepsilon/2 - \alpha (4 d - 2)} + |B_{2 N}|\, e^{- N^{\epsilon/d}}
\end{equation*}
which implies \eqref{SNOc}.
\end{proof}

In conclusion, in view of \twoeqref{eta+}{S1}, \treeeqref{n0}{tauN}{Nt} and \twoeqref{SNO}{SNOc}, we obtain that 
\begin{equation*}
\PP \left(P^o_\omega (X_{2 n} = o) > 1/m\right) \lesssim m\, \max \{n^{\varepsilon/2 - 2 - \alpha (4 d - 2)}, n^{- d/2}\} + o (1) 
\end{equation*}
The demonstration is done.


\section{\bf The case \eqref{UP}}\label{}

\begin{proof}[Proof of the case $\alpha > 0$]
For this first general case, we don't need to assume \eqref{UP}, as the result remains unconditionally true for conductances in $[1, \infty)$.
Let us start by adapting the notation of the previous model with the new one.
Choose $\theta > 1$ such that conductances less than $\theta$ percolate and let $\scrC^\theta$ denote the unique infinite connected cluster of conductances in $[1, \theta]$ which we can suppose to contain $o$ without loss of generality.
Set $A_\theta (t) = \int^t_0 \1_{\{X_s \in \scrC^\theta\}}\, {\rm d} s$
and let $A^{- 1}_\theta (t)$ be its right-continuous inverse. 
Then, define $X^\theta_{t} = X_{A^{- 1}_\theta (t)}$, the so-called trace of $X$ on $\scrC^\theta$. 
This is a random walk that has normal transition probabilities. 
Indeed, for a.e. $\omega$, for any $x, y \in \scrC^\theta$, we have (cf. \cite{BC}),
\begin{equation}
\label{Zttp}
P^x_\omega (X^\theta_t = y)  \lesssim t^{- d/2}
\end{equation}
Additionally, we have by ergodicity \cite{ABDH} that $\PP-$a.s., 
\begin{equation}
\label{Atht}
t^{- 1} A_\theta (t) \xrightarrow[t \to \infty]{} \PP (o \in \scrC_\theta), \quad P^ o_\omega-a.s.
\end{equation}

Now, observe that
\begin{equation}
\label{fs-ineq}
P^o_\omega(X_t = o) \leq \frac{2}{t}\int^{t}_{t/2} P^o_\omega(X_v = o)\, \text d v = \frac{2}{t}\, E^o_\omega\left(\int^{t}_{t/2}\1_{\{X_v = o\}}\, \text d v\right)
\end{equation}
where we used that the return probability is decreasing. 
Since $A_\theta (t)$ is a positive continuous increasing function of the time, we obtain by operating a time change in the last integral in \eqref{fs-ineq},
$$
\int^{t}_{t/2}\1_{\{X_v = o\}}\text d v
= \int^{A_\theta (t)}_{A_\theta (t/2)}\1_{\{X^\theta_u = o\}}\text du
\le \int^{\infty}_{A_\theta (t/2)}\1_{\{X^\theta_u = o\}}\text du
$$
In addition, by Fubini theorem and \eqref{Zttp}, we have, for any $s > 0$ and for $d \ge 3$, 
$$
E^o_\omega\left(\int^\infty_s \1_{\{X^\theta_u = o\}}\, \text d u\right) \le 
\int^\infty_s  u^{- d/2}\, {\rm d} u < \infty
$$
which implies that 
$$
\int^\infty_s  \1_{\{X^\theta_u = o\}}\, {\rm d} u \xrightarrow[s \to \infty]{} 0, \quad P^ o_\omega-a.s.
$$

Thus, it comes by monotonicity and \eqref{Atht} that
$$
E^o_\omega\left(\int^\infty_{A_\theta (t/2)} \1_{\{X^\theta_u = o\}}\, \text d u\right) \eqqcolon \epsilon_t 
\xrightarrow[t \to \infty]{} 0 
$$
Hence, we get that $P^o_\omega(X_t = o)\lesssim \epsilon_t\, t^{- 1}$.
\end{proof}

\begin{proof}[Proof of the case $\alpha > 1$]
For this integrable conductance case, we have on the one hand an  immediate normal lower bound that stems from that of $Y$. Indeed, for $d \ge 2$, we have $P^o_\omega (Y_t = o) \simeq t^{- d/2}$ (cf. \cite{ABDH}). Besides, $A (t) \ge t$ because of $\pi (x) \ge 2 d$, which yields $A^{- 1} (t) \le t$. 
Since the return probability is decreasing, we obtain that
$$
P^ o_\omega (Y_{A^{- 1} (t)} = o) \ge \inf_{s \ge t} P_\omega^o \big( Y_{s} = o\big) \gtrsim t^{- d/2}
$$

On the other hand, let us deal with the upper bound. 
Suppose that $\mu = \E (\pi (o)) < \infty$, or $\alpha > 1$ under the condition \eqref{UP}.
Then, by the ergodic theorem (cf. \cite{ABDH}), we have $\PP-$a.s., 
\begin{equation}
\label{Atpi}
t^{- 1} A (t) \xrightarrow[t \to \infty]{} \E (\pi (o)) \eqqcolon \mu, \quad P^{o}_\omega-a.s.
\end{equation}
Thereby, we intuitively expect $X_t$ to behave like $Y_{t/\mu}$ since $X_t = Y_{A^{- 1} (t)}$ and therefore $X$ would have normal return probabilities. 

Formally, for any $c > 0$, observe that
\begin{eqnarray}
P_\omega^o \big( X_t = o \big)
&=&
P_\omega^o \big(Y_{A^{-1} (t)} = o, A^{-1} (t) \ge c\,t\big)
\nonumber
\\
&&\qquad\qquad
+ P_\omega^o \big( X_t = o, A^{-1} (t) < c\, t\big)
\nonumber
\\
&\le&
\sup_{s \ge c\, t} P_\omega^o \big( Y_{s} = o\big) + P_\omega^o \big( X_t = o, A^{-1} (t) < c\, t\big)
\nonumber
\\
\label{rev}
&\lesssim& 
t^{- d/2} + P_\omega^o \big( X_t = o, A^{-1} (t) < c\, t\big)
\end{eqnarray}
So, it remains to estimate the second term on the right-hand side of \eqref{rev}. Note that $\mu \ge 2 d$ because of $\pi (x) \ge 2 d$ and choose 
\begin{equation}
\label{c}
c \in \left(\frac{1}{\mu (\mu + 1)},\frac{1}{\mu + 1}\right)
\end{equation}
Furthermore, set $\nu = 1 + {1}/{\mu}, c^* = c \,\nu, t^* = {t}/\mu$
and remark that $t^* < t/(2 d)$. Since $A (t)$ is (increasingly) convex in $t$ with $A (t) > 2 d t$, we have that $A^{- 1} (t)$ is (increasingly) concave in $t$ with $A^{- 1} (t) < t/(2 d)$, which easily yields
\begin{equation*}
A^{- 1} (t) < c \,t \Longrightarrow A^{- 1} (\nu t) \le c^* t
\end{equation*}
Consequently, we obtain 
\begin{eqnarray}
\label{A-t0}
P_\omega^o \big( X_t = o, A^{-1} (t) < c\, t\big) 
\le 
P_\omega^o \big( X_t = o\big) - P_\omega^o \big( X_t = o, A^{- 1} (\nu t) > c^* t\big)
\end{eqnarray}
But, we also have
$$
A^{-1} (t^*) > c^* t \Longrightarrow A^{- 1} (\nu t) > c^* t 
$$
which, by the Markov property, implies that
$$
P_\omega^o \big( X_t = o, A^{- 1} (\nu t) > c^* t\big) 
\ge P^ o_\omega (X_t = o)   P^ o_\omega \left(A^{-1} (t^*) > c^* t\right)
$$
Then, by using this last inequality in \eqref{A-t0}, it comes that the second term on the right-hand side of \eqref{rev} is less than
$$
P_\omega^o \big( X_t = o\big) \left(1 - P^ o_\omega \left(A^{-1} (t^*) > (\mu c^*)\, t^*\right)\right)
$$
Besides, by symmetry, we have
\begin{equation}
\label{A-t}
A^{-1} (t^*)/t^* > \mu c^* \Longleftrightarrow A (t^*)/t^* < 1/(\mu c^*)
\end{equation}
But, $1/(\mu c^*) > \mu$ because of \eqref{c}. Accordingly, \eqref{Atpi} gives that
$$
P^ o_\omega \left(A^{-1} (t^*) > (\mu c^*)\, t^*\right) \xrightarrow[t \to \infty]{} 1
$$

Thus, by substitution in \eqref{rev}, the desired standard upper bound comes from
$$
P_\omega^o \big(X_t = o\big) \lesssim t^{- d/2} + P_\omega^o \big(X_t = o\big)\, o (1)
$$  
\end{proof}


\medskip
\begin{proof}[Proof of the case $\alpha < 1/2$]
The strategy is the same as for  the case \eqref{LP}, however,  the calculations require an adaptation and here is a sketch of the proof.

First, let $m = [\Theta_t^{- 1} t^{\varepsilon}] \simeq t^{\alpha + \varepsilon}$ with an arbitrary $\varepsilon > 0$. Then, as in \eqref{Lm}, we obtain for $t$ large enough,
\begin{equation*}
\label{Lamn}
P^o_\omega \left(\sum^m_{k = 1} \1_{V_k} \ge 1\right) \ge 1 - e^{- t^{\varepsilon}}
\end{equation*}
where $V_k$ is a trapping configuration such as $\pi (X_{H_k}) \simeq 2 d \theta$ with an adjacent trap outside $B_k$ as defined above.

Second, if $\alpha + \varepsilon < 1/2$, we claim that
\begin{equation}
\label{Hmt}
\lim_{t \to \infty} P^o_\omega (H_m \ge t) = 0
\end{equation}
Indeed, if $H^\theta_m$ denotes the hitting time of $\partial B_m$ by $X^\theta$, a standard argument (cf. \cite[Lemma~2.1]{BKM}) which uses \eqref{Zttp} and the Markov property implies that 
\begin{equation}
\label{EHkth}
E^ o_\omega \left(H^\theta_m\right) \lesssim m^2
\end{equation} 
Then, the claim follows from
\begin{equation*}
\{H_m \ge t \} \subset \{H_m^\theta \ge A_\theta (t)\} \subset  \{H_m^\theta \ge t/2\} 
\cup \{A_\theta (t) \le t/2\}
\end{equation*}
which gives
$$
P^ o_\omega \big(H_m \ge  t\big) 
\le
P^ o_\omega \big(H^\theta_{m} \ge t/2\big) 
+ P^ o_\omega \big(A_\theta (t) \le t/2\big)
\lesssim
t^{2 (\alpha + \varepsilon) - 1} + o (1)
$$
where we used Markov inequality with \eqref{EHkth} and \eqref{Atht} with $\theta$ large enough.

Next, remark that the probability for the random walk to step outside a trap of depth $t$ is less than $\theta/(2 d t)$ ($\theta\ll t$). Thus, as in Lemma~\ref{XinS} and thanks to \eqref{Hmt} with $\alpha + \varepsilon < 1/2$, we have for $t$ large enough, for each $k = 1, \ldots, m$, 
\begin{equation*}
\label{stayintrap+}
P^ o_\omega \big(X_t \in \mathbf e_k \mid V_{k}\big) \ge \frac{p_t}{2 d}\, \left(1 - \frac{\theta}{2 d t}\right)^{t} 
\gtrsim \frac{e^{- \theta/(2 d)}}{2 d}
\end{equation*}
where $p_t$ is the probability to have less than $t$ jumps before time $t$; this Poisson distribution probability tends to $1$.
For the rest of the proof, we can mimic the steps of the \mbox{case \eqref{LP}}. 
\end{proof}

\vspace{-0,3cm}
\section{\bf Acknowledgments}

The author would like to thank Pierre, Marek, Takashi and Nina.



\begin{thebibliography}{00}



\bibitem{ABDH}
S.~Andres, M.T.~Barlow, J.-D.~Deuschel and B.M.~Hambly (2013).
\textit{Invariance principle for the random conductance model}, 
Probab. Theory Rel. Fields., {\bf 156} (2013) no.~3, 535--580.


\bibitem{BC}
M. T. Barlow\ and\ J. \v Cern\'y, \textit{Convergence to fractional kinetics for random walks associated with unbounded conductances}, Probab. Theory Related Fields {\bf 149} (2011), no.~3-4, 639--673.


\bibitem{BD} 
M. T. Barlow\ and\ J.-D. Deuschel, \textit{Invariance principle for the random conductance model with unbounded conductances}, Ann. Probab. {\bf 38} (2010), no.~1, 234--276.


\bibitem
{BBHK}
N. Berger,\ M. Biskup, C. E. Hoffman, G. Kozma, \textit{Anomalous heat kernel decay for random walk among bounded random conductances}, Ann. Inst. Henri Poincar\'e Probab. Stat. {\bf 44} (2008), no.~2, 374--392.



\bibitem
{M-review}
M. Biskup  (2011). 
\textit{Recent progress on the Random Conductance Model}. 
\textit{Prob. Surveys} \textbf{8} 294--373.

%
\bibitem
{BiBo}
M. Biskup\ and\ O. Boukhadra, \textit{Subdiffusive heat kernel decay in four-dimensional i.i.d. random conductance models}, 	J. London Math. Soc. {\bf 86} (2012), no.~2, 455--481.


\bibitem
{UCLA-team}
M. Biskup\ et al., \textit{Trapping in the random conductance model}, J. Stat. Phys. {\bf 150} (2013), no.~1, 66--87. 



\bibitem
{B1}
O. Boukhadra, \textit{Heat kernel estimates for random walk among random conductances with heavy tail}, Stochastic Process. Appl. {\bf 120} (2010), no.~2, 23--27.


\bibitem
{B-SPL}
O. Boukhadra, \textit{On heat kernel decay for random conductance model}, Statistics and Probability Letters, {\bf 133} (2018)


\bibitem{B2}
O. Boukhadra, Standard spectral dimension for the polynomial lower tail random conductances model, Electron. J. Probab. {\bf 15} (2010), no. 68, 2069--2086.

%
\bibitem
{BKM}
O. Boukhadra\, T. Kumagai and P. Mathieu, \textit{Local CLT for the polynomial lower tail RCM}, J. Math. Soc. Japan Vol. 67, No. 4 (2015) pp. 1413–1448. 

\bibitem{BP}
M. Biskup\ and\ T. M. Prescott, \textit{Functional CLT for random walk among bounded random conductances}, 
{\em Electron. J. Probab. \bf 12} (2007), no. 49, 1323--1348.

%
\bibitem{C}
J. \v Cern\'y, \textit{On two-dimensional random walk among heavy-tailed conductances}, Electron. J. Probab. {\bf 16} (2011), no. 10, 293--313.
%

\bibitem
{FM}
L. R. G. Fontes\ and\ P. Mathieu, \textit{On symmetric random walks with random conductances on ${\Bbb Z}\sp d$}, Probab. Theory Related Fields {\bf 134} (2006), no.~4, 565--602.
%


\bibitem
{G}
G. Grimmett, {\it Percolation}, second edition, Grundlehren der Mathematischen Wissenschaften, 321, Springer, Berlin, 1999. 


\bibitem
{K}
T. Kumagai, \emph{Random Walks on Disordered Media and their Scaling Limits}, 
Lect. Notes in Math., {\bf 2101}, \'Ecole d'\'Et\'e de Probabilit\'es de Saint-Flour XL--2010. 
Springer, New York, (2014). 


\bibitem
{MR}
P. Mathieu\ and\ E. Remy, \textit{Isoperimetry and heat kernel decay on percolation clusters}, 
{\em Ann. Probab. \bf 32} (2004), no.~1A, 100--128. 


\end{thebibliography}
\end{document}